\documentclass[11pt]{article}
\usepackage[table]{xcolor}
\usepackage{graphicx}
\usepackage{amsmath,amsfonts,amssymb}
\usepackage{graphicx}
\usepackage{multirow}
\usepackage{tikz,pgf}
\usepackage{url}
\usepackage{amsmath}
\usepackage{graphicx}
\usepackage{epsfig}
\usepackage{epstopdf}
\usepackage{color}
\usepackage{enumitem}
\usepackage{cite}
\def\disp{\displaystyle}

\def\tto{\;{\lower 1pt \hbox{$\rightarrow$}}\kern -10pt
\hbox{\raise 2pt \hbox{$\rightarrow$}}\;}
\def\Hat{\widehat}

\def\Bar{\overline}
\def\ra{\rangle}
\def\la{\langle}

\def\epsilon{\varepsilon}
\def\B{\Bbb B}
\def\h{\hfill\Box}
\def\R{\Bbb R}

\def\ox{\bar{x}}

\def\oy{\bar{y}}
\def\oz{\bar{z}}

\def\co{\mbox{\rm co}}

\def\ri{\mbox{\rm ri}}

\def\gph{\mbox{\rm gph}}
\def\aff{\mbox{\rm aff}}
\def\epi{\mbox{\rm epi}}

\def\dom{\mbox{\rm dom}}
\def\aff{\mbox{\rm aff}}

\def\cone{\mbox{\rm cone}}
\def\rge{\mbox{\rm rge}}

\def\h{\hfill\square}

\def\ph{\varphi}
\def\emp{\emptyset}

\def\oR{\Bar{\R}}

\def\ph{\varphi}
\def\emp{\emptyset}

\def\oR{\Bar{\R}}

\setlist[enumerate,1]{itemsep=0.0ex,parsep=0.5ex,label={\rm(\alph*)},leftmargin=*, align=left}
\newcounter{lk}

\usepackage[colorlinks=true]{hyperref}

\begin{document}
\begin{center}
{\sc\bf Revisiting Rockafellar's Theorem on Relative Interiors of Convex Graphs\\ with Applications to Convex Generalized Differentiation}\\[1ex]
{\sc Dang Van Cuong} \footnote{Department of Mathematics, Faculty of Natural Sciences, Duy Tan University, Da Nang, Vietnam (dvcuong@duytan.edu.vn). This research is funded by Vietnam National Foundation for Science and Technology Development (NAFOSTED) under grant number 101.02-2020.20.}, {\sc B. S. Mordukhovich}\footnote{Department of Mathematics, Wayne State University, Detroit, Michigan 48202, USA (boris@math.wayne.edu). Research of this author was partly supported by the USA National Science Foundation under grant DMS-1808978, by the Australian Research Council under grant DP-190100555, and by Project 111 of China under grant D21024.}, {\sc Nguyen Mau  Nam}\footnote{Fariborz Maseeh Department of Mathematics and Statistics, Portland State University, Portland, OR
97207, USA (mnn3@pdx.edu). Research of this author was partly supported by the USA National Science Foundation under grant DMS-2136228.},
 {\sc G. Sandine}\footnote{Fariborz Maseeh Department of
 Mathematics and Statistics, Portland State University, Portland, OR 97207, USA (gsandine@pdx.edu).}\\[2ex]
  {\bf Dedicated to Roger Wets with the highest respect}
\end{center}
\small{\bf Abstract.} In this paper we revisit a theorem by Rockafellar on representing the relative interior of the graph of a convex set-valued mapping in terms of the relative interior of its domain and function values. Then we apply this theorem to provide a simple way to prove many calculus rules of generalized differentiation for set-valued mappings and nonsmooth functions in finite dimensions. Using this important theorem by Rockafellar allows us to improve some results on generalized differentiation of set-valued mappings in \cite{bmncal} by replacing the relative interior qualifications on graphs with qualifications on domains and/or ranges.\\[1ex]

\newtheorem{Theorem}{Theorem}[section]
\newtheorem{Proposition}[Theorem]{Proposition}
\newtheorem{Remark}[Theorem]{Remark}
\newtheorem{Lemma}[Theorem]{Lemma}
\newtheorem{Corollary}[Theorem]{Corollary}
\newtheorem{Definition}[Theorem]{Definition}
\newtheorem{Example}[Theorem]{Example}
\renewcommand{\theequation}{\thesection.\arabic{equation}}
\normalsize

\section{Introduction}
\setcounter{equation}{0}

The notion of {\em relative interior} for convex sets in finite-dimensional spaces goes back to Steinitz \cite{s} and then has been systematically studied and applied in finite-dimensional convex analysis and related areas; see, e.g., the seminal monograph by Rockafellar \cite{r} and the subsequent publications including the books by Borwein and Lewis \cite{Borwein2000} and by Hiriart-Urruty and Lemar\'echal \cite{HU1} with further references and commentaries. In contrast to the interior, the relative interior is {\em nonempty} for any nonempty convex set in $\R^n$, while the latter notion shares many important properties of the interior being very useful in applications. Relative interiors play a crucial role in many aspects of  convex analysis and optimization in finite dimensions such as convex separation, generalized differential calculus, Fenchel conjugate, and Fenchel and Lagrange duality; see, e.g., \cite{Bauschke2011,Bertsekas2003,Borwein2000,Boyd2004,g,HU1,KK2013,bmn,mn22,pr,rw} and the references therein.

Among many important results involving relative interiors is the theorem by Rockafellar \cite[Theorem~6.8]{r} allowing us to represent the relative interior of a convex set $G$ in $\R^n\times \R^m$ in terms of the the relative interiors of the image $D$ of $G$ under the projection mapping $(x, u)\to x$ and  the set $S(x):=\{u\in \R^m\; |\; (x, u)\in G\}$ for $x\in D$. In the language of set-valued analysis, this theorem gives a representation of the relative interior of the graph of a convex set-valued mapping in terms of the relative interiors of its domain and the function's values. The first proof of this result was provided in the aforementioned book by Rockafellar as a consequence of several other results involving relative interiors of convex sets. A more self-contained proof  was given by Rockafellar and Wets \cite[Proposition~2.43]{rw}. In a recent paper \cite{CBN}, we provided the third proof of this important theorem and explored its generalization to locally convex topological vector spaces using a generalized relative interior concept called the \emph{quasi-relative interior} introduced in \cite{bl}.

The first goal of the present paper is to revisit Rockafellar's theorem and derive a new result on relative interiors of graphs of \emph{generalized epigraphical mappings}. Then we employ these results and the geometric approach to convex analysis developed in \cite{bmncal,mn22} to provide a simple way to access many calculus rules of generalized differentiation for set-valued mappings under new qualification conditions. The usage of Rockafellar's theorem and related developments allows us to improve, in particular, a number of calculus rules obtained in \cite{bmncal} for coderivatives of convex set-valued mappings. Our developments have a great potential for further implementations in the field of set-valued optimization; see, e.g., the books \cite{AAK,m18} and the references therein for this and related areas of optimization theory and applications.

This paper is organized as follows. Section~2 contains basic concepts of convex analysis in finite dimensions used throughout the paper. In Section~3, we revisit Rockafellar's theorem on relative interiors of convex graphs and derive a number of new results with detailed proofs. Section~4 is devoted to employing this theorem and the geometric approach to convex analysis in the study of generalized differentiation for convex set-valued mappings and nonsmooth functions in finite dimensions. Some applications to convex generalized equations, convex constraint systems, and optimal value functions are presented in Section~5. In Section~6, we develop the convex coderivative calculus for set-valued mappings obtained under relative interior qualification conditions imposed on domains or ranges. These developments significantly improve calculus rules in \cite[Section~11]{bmncal} under relative interior qualification conditions imposed on graphs. Throughout the paper, we use standard notation of convex analysis in finite dimensions; see \cite{Borwein2000,HU1,mn22,r}. In particular,  $\la x, y\ra$ denotes the inner product of $x, y\in\R^n$;  $\B(x;\gamma)$ signifies the closed ball centered at $x$ with radius $\gamma\geq 0$; the  closure of a set $\Omega\subset\R^n$ is denoted by $\Bar{\Omega}$; the convex hull of a set $\Omega$ is $\co(\Omega)$; the cone generated by a set $\Omega$ is $\cone(\Omega):=\{tw\; |\; t\geq 0, w\in \Omega\}$.

\section{Preliminaries}
\setcounter{equation}{0}
In this section, we recall a number of concepts and results of convex analysis in finite dimensions used throughout the paper; see, e.g., \cite{Borwein2000,HU1,mn22,r} and the references therein.

A subset $\Omega$ of $\R^n$ is called {\em convex} if
\begin{equation*}
\lambda x+(1-\lambda)y\in \Omega\; \mbox{\rm  for all }x,y\in\Omega\; \mbox{\rm  and }\lambda\in(0,1).
\end{equation*}
It follows directly from the definition that $\Omega$ is convex if and only if for any two points $x,y\in\Omega$, the line segment $(x, y):=\{\lambda x+(1-\lambda)y\;|\;\lambda\in (0,1)\}$ is a subset of $\Omega$. A subset $\Omega$ of $\R^n$ is called  {\em affine} if for any $x,y\in \Omega$ and for any $\lambda\in\R$ we have
\begin{equation*}
\lambda x+(1-\lambda)y\in \Omega,
\end{equation*}
which means that $\Omega$ is affine if and only if the line through any two points $x,y\in \Omega$ is a subset of $\Omega$. It follows directly from the definition that any affine set is a convex set. In addition,  the intersection of any collection of affine sets is an affine set and thus allows us to define the {\em affine hull} of a set $S$ by
\begin{equation*}
\aff(S):=\bigcap\{\Omega\; |\;\Omega\;\text{ is affine and }\;S\subset \Omega\}.
\end{equation*}
The {\em relative interior} $\ri(\Omega)$ of a set $\Omega$ in $\R^n$  is defined as its interior within the affine hull of $\Omega$, i.e., by
\begin{equation*}
\ri(\Omega):=\{x\in \Omega\ |\ \exists\, \gamma>0\ \text{satisfying}\ \B(x;\gamma)\cap\aff(\Omega)\subset\Omega\}.
\end{equation*}

Let $\Omega$ be a nonempty convex subset of $\R^n$ with $\ox\in\Omega$. The {\em normal cone} to $\Omega$ at $\ox$ is
\begin{equation*}\label{nc}
N(\ox;\Omega):=\big\{v\in \R^n\;\big|\;\la v,x-\ox\ra\le 0\;\text{ for all }\;x\in\Omega\big\}
\end{equation*}
with $N(\ox;\Omega):=\emp$ if $\ox\notin\Omega$.

The following theorem provides several characterizations for the relative interior of a convex set; see, e.g., \cite[Theorem~2.2]{CBN}.
\begin{Theorem}{\bf(characterizations of relative interior for convex sets in $\R^n$).}\label{pts2} Let $\Omega$ be a nonempty convex set in $\R^n$ and let $\ox\in \R^n$. The following properties are equivalent:
\begin{enumerate}
\item $\ox\in\ri(\Omega)$.
\item  $\ox\in\Omega$ and for every $x\in\Omega$ with $x\ne\ox$ there exists $u\in\Omega$ such that
$\ox\in(x,u)$.
\item  $\ox\in\Omega$ and $\cone(\Omega-\ox)$ is a linear subspace of $\R^n$.
\item  $\ox\in\Omega$ and $\overline{\cone}(\Omega-\ox)$ is a linear subspace of $\R^n$.
\item  $\ox\in\Omega$ and the normal cone $N(\ox;\Omega)$ is a subspace of $\R^n$.
\end{enumerate}
\end{Theorem}
The relative interior possesses several nice algebraic and topological properties, some of which are presented in the theorem below; see, e.g., \cite{r}.
\begin{Theorem}{\bf (properties of relative interiors).} \label{Ri_nonem}
Let $\Omega$ and $\Omega_i$ for $i=1, \ldots, m$ be nonempty convex subsets of $\R^n$. Then
\begin{enumerate}
\item $\ri(\Omega)$ is nonempty and convex.
\item $[a,b)\subset\ri(\Omega)$ for any $a\in\ri(\Omega)$ and $b\in\Bar{\Omega}$, where $[a,b):=\{ta+(1-t)b\;|\;0<t\le 1\}$ defines the half-open interval connecting $a, b\in \R^n$.
    \item $\Bar{\Omega}=\Bar{\ri(\Omega)}$ and $\ri(\Bar{\Omega})=\ri(\Omega)$.
    \item $\ri(\ri(\Omega))=\ri(\Omega)$.
    \item $\ri(\sum_{i=1}^m\Omega_i)=\sum_{i=1}^m \ri(\Omega_i)$.
    \item $\ri (A(\Omega))=A(\ri(\Omega))$, where $A\colon \R^n\to \R^m$ is a linear mapping.
    \item $\ri(\cap_{i=1}^m\Omega_i)=\cap_{i=1}^m \ri(\Omega_i)$ provided that $\cap_{i=1}^m \ri(\Omega_i)\neq\emptyset$.
\end{enumerate}
\end{Theorem}
It is worth noting that the relative interior may not inherit all properties of the interior. For example, for two nonempty convex sets $\Omega_1$ and $\Omega_2$ in $\R^n$ with $\Omega_1\subset \Omega_2$, it is not true in general that $\ri(\Omega_1)\subset \ri(\Omega_2)$.

Another important role of the relative interior is in the study of convex proper separation in $\R^n$. Recall  that two nonempty convex sets $\Omega_1,\Omega_2\subset \R^n$ can be {\em properly separated} if there exists $v\in \R^n$ for which the following two inequalities hold:
\begin{equation}\label{ps1}
\sup\big\{\la v,w_1\ra\;\big|\;w_1\in\Omega_1\big\}\le\inf\big\{\la v,w_2\ra\;\big|\;w_2\in\Omega_2\big\},
\end{equation}
\begin{equation}\label{ps2}
\inf\big\{\la v,w_1\ra\;\big|\;w_1\in\Omega_1\big\}<\sup\big\{\la v,w_2\ra\;\big|\;w_2\in\Omega_2\big\}.
\end{equation}
Observe that condition \eqref{ps1} can be equivalently rewritten as
\begin{equation*}
\la v,w_1\ra\le\la v,w_2\ra\;\mbox{\rm whenever }w_1\in\Omega_1,\;w_2\in\Omega_2,
\end{equation*}
while \eqref{ps2} means that there exist $\Hat{w}_1\in \Omega_1$ and $\Hat{w}_2\in \Omega_2$ such that
\begin{equation*}
\la v,\Hat{w}_1\ra<\la v,\Hat{w}_2\ra.
\end{equation*}
As a central theorem of convex analysis in finite dimensions, the following theorem uses the relative interior to provide necessary and sufficient conditions for properly separating two nonempty convex sets; see, e.g., \cite[Theorem~11.3]{r}.

\begin{Theorem}{\bf(relative interior and proper separation in finite dimensions).}\label{pst1} Let $\Omega_1$ and $\Omega_2$ be two nonempty convex subsets of $\R^n.$ Then $\Omega_1$ and $\Omega_2$ can be properly separated if and only if $\ri(\Omega_1)\cap \ri(\Omega_2)=\emptyset.$
\end{Theorem}
The relative interior  plays a crucial role in many other issues of convex analysis. For instance, a direct application of Theorem~\ref{pst1} for $\Omega$ and the single-point set $\{\ox\}$ shows that if $\ox\in \Omega\setminus \ri(\Omega)$, then $N(\ox; \Omega)\neq\{0\}$. The relative interior and the proper separation theorem can also be used in the statement and proof of the normal cone intersection rule in the theorem below; see, e.g., \cite[Corollary 23.8.1]{r} for more details.

\begin{Theorem}{\bf (normal cone intersection rule in finite dimensions).}\label{ncr} Let $\Omega_1,\ldots,\Omega_m\subset\R^n$ be convex sets satisfying the relative interior condition
\begin{eqnarray*}\label{qc1}
\bigcap_{i=1}^m\ri(\Omega_i)\ne\emp,
\end{eqnarray*}
where $m\ge 2$. Then we have the normal cone intersection rule
\begin{equation*}\label{sum2}
N\Big(\bar{x};\bigcap_{i=1}^m\Omega_i\Big)=\disp\sum_{i=1}^m N(\bar{x};\Omega_i)\;\mbox{\rm for all }\;\bar{x}\in\bigcap_{i=1}^m\Omega_i.
\end{equation*}
\end{Theorem}
Given a set-valued mapping $F\colon \R^n\tto \R^m$, the  \emph{graph} of $F$ is the set
\begin{equation*}
\gph(F):=\big\{(x,y)\in \R^n\times \R^m\;\big|\;y\in F(x)\big\},
\end{equation*}
and it is called {\em convex} if $\gph(F)$ is a convex subset of the product space $\R^n\times \R^m$. We also consider the {\em domain} and {\em range} of $F$ defined by
\begin{equation*}
\dom(F):=\big\{x\in \R^n\;\big|\;F(x)\ne\emp\big\}\; \mbox{ and }\;\rge(F):=\bigcup_{x\in \R^n}F(x),
\end{equation*}
respectively. It is easy to see that if $F$ is a convex set-valued mapping, then $\dom(F)$ and $\rge(F)$ are convex sets as well.

The \emph{coderivative} of a convex set-valued mapping $F$ at $(\ox, \oy)\in \gph(F)$ is defined by
\begin{equation}\label{cod}
D^*F(\ox,\oy)(v):=\big\{u\in\R^n\;\big|\;(u,-v)\in N((\ox,\oy);\gph(F))\big\},\; v\in\R^m.
\end{equation}
The coderivative can be used to define the  \emph{subdifferential} of an extended-real-valued convex function $f\colon \R^n \to (-\infty, \infty]$. Given $\ox\in \dom(f):=\{x\in \R^n\; |\; f(x)<\infty\}$, define
\begin{equation}\label{cvsd}
\partial f(\ox):=D^*E_f(\ox, f(\ox))(1)=\big\{v\in \R^n\; \big|\; \la v, x-\ox\ra\leq f(x)-f(\ox)\; \mbox{\rm for all }x\in \R^n\big\},
\end{equation}
where $E_f(x):=[f(x), \infty)$ for $x\in \R^n$ is the \emph{epigraphical mapping/multifunction} associated with $f$ with $\gph(E_f)=\epi(f):=\{(x, \lambda)\in \R^n\times \R\; |\; f(x)\leq \lambda\}$.

\section{Rockafellar's Theorem on Relative Interiors of Convex Graphs}
\setcounter{equation}{0}

In this section, we revisit the aforementioned theorem by Rockafellar on representing relative interiors of graphs of convex set-valued mappings.

\begin{Theorem}{\bf(Rockafellar's theorem on relative interiors of convex graphs).}\label{TheoRoc1} Let $F\colon\R^n\tto\R^m$ be a convex set-valued mapping. Then we have the representation
\begin{equation}\label{roc}
\ri\big(\gph(F)\big)=\big\{(x,y)\in\R^n\times\R^m\;\big|\;x\in\ri\big(\dom(F)\big),\;y\in\ri\big(F(x)\big)\big\}.
\end{equation}
\end{Theorem}
{\bf Proof.} We first prove the inclusion ``$\subset$" in \eqref{roc}. Consider the projection mapping $\mathcal{P}\colon \R^n\times \R^m\to \R^n$ given by
\begin{equation*}
\mathcal{P}(x, y):=x\; \mbox{\rm for }(x, y)\in \R^n\times \R^m.
\end{equation*}
It follows from Theorem \ref{Ri_nonem}(f) that
\begin{equation}\label{proj}
\mathcal{P}(\ri(\gph(F))=\ri(\mathcal{P}(\gph(F)))=\ri(\dom(F)).
\end{equation}
Now, take any  $(\ox,\oy)\in\ri(\gph(F))$ and get from \eqref{proj} that   $\ox\in\ri(\dom (F))$. Since $(\ox,\oy)\in\ri(\gph(F))\subset \gph(F)$, we have $\oy\in F(\ox)$. Fix any $y\in F(\ox)$ with $y\neq \oy$. Then $(\ox, y)\in \gph(F)$ with $(\ox, y)\neq (\ox, \oy)$. By the equivalence of (a) and (b) from Theorem \ref{pts2}, there exists $(u,z)\in \gph(F)$ and $t\in (0, 1)$ such that
\begin{equation*}
(\ox, \oy)=t(\ox, y)+(1-t)(u, z).
\end{equation*}
Then $\ox=t\ox+(1-t)u$, which implies $(1-t)\ox=(1-t)u$ and so $\ox=u$. In addition, $\oy=ty+(1-t)z\in (y, z)$, where $z\in F(\ox)$. Using the equivalence of (a) and (b) from Theorem \ref{pts2} again yields $\oy \in\ri(F(\ox))$.

To verify next the reverse inclusion in \eqref{roc}, fix $\ox\in\ri(\dom(F))$ and $\oy\in\ri(F(\ox))$. Arguing by contradiction, suppose that $(\ox,\oy)\notin\ri(\gph(F))$ and then find by Theorem~\ref{pst1} a pair $(u,v)\in\R^n\times\R^m$ such that
\begin{equation}\label{1}
\la u,x\ra+\la v,y\ra\le\la u,\ox\ra+\la v,\oy\ra\;\mbox{ whenever }\;(x,y)\in\gph(F).
\end{equation}
In addition, it follows from the proper separation of $\{(\ox,\oy)\}$ and $\gph(F)$ that there exists a pair $(x_0,y_0)\in\gph(F)$ satisfying
\begin{equation}\label{2}
\la u,x_0\ra+\la v,y_0\ra<\la u,\ox\ra+\la v,\oy\ra.
\end{equation}
Letting $x:=\ox$ in \eqref{1} yields $\la v,y\ra\le\la v,\oy\ra$ for all $y\in F(\ox)$.  Since $\ox\in\ri(\dom(F))$ and $x_0\in \dom(F)$, we deduce from Theorem~\ref{pts2} that there exists $\tilde{x}\in\dom(F)$ such that $\ox=tx_0+(1-t)\tilde{x}$ for some $t\in(0,1)$, which is is true even if $x_0=\bar{x}$. Choose $\tilde{y}\in F(\tilde{x})$ and consider the convex combination
\begin{equation*}
y^\prime:=ty_0+(1-t)\tilde{y},
\end{equation*}
where $y^\prime\in F(\ox)$ since $\gph(F)$ is convex. Since $(\tilde{x}, \tilde{y})\in\gph(F)$ we use \eqref{1} and \eqref{2} to get
\begin{align*}
&\la u,\tilde{x}\ra+\la v,\tilde{y}\ra\le\la u,\ox\ra+\la v,\oy\ra,\\
&\la u,x_0\ra+\la v,y_0\ra<\la u,\ox\ra+\la v,\oy\ra.
\end{align*}
Multiplying the first inequality above by $1-t$ and the second one by $t$, and then adding them together give us the condition
\begin{equation*}
\la u,\ox\ra+\la v,y^\prime\ra<\la u,\ox\ra+\la v,\oy\ra,
\end{equation*}
which yields $\la v,y^\prime\ra<\la v,\oy\ra$. From this along with~\eqref{1} when $x=\ox$, we conclude that the sets $\{\oy\}$ and $F(\ox)$ can be properly separated. Applying Theorem~\ref{pst1} tells us that $\oy\notin\ri(F(\ox))$, a contradiction that verifies $(\ox,\oy)\in\ri(\gph(F))$. $\h$

Given a set-valued mapping $F\colon \R^n\tto \R^m$, recall that the {\em inverse} of $F$ is the set-valued mapping $F^{-1}\colon \R^m\tto \R^n$ defined by
\begin{equation*}
F^{-1}(y):=\big\{x\in \R^n\; \big|\; y\in F(x)\big\}, \; y\in \R^m.
\end{equation*}
The next corollary is a direct consequence of Theorem \ref{TheoRoc1}.
\begin{Corollary}\label{rirange}{\bf (relative interiors of convex graphs and ranges).} Let $F\colon \R^n\tto \R^m$ be a convex set-valued mapping. If $(\ox, \oy)\in \ri(\gph(F))$, then $\oy\in \ri(\rge(F))$.
\end{Corollary}
{\bf Proof.} It is not hard to show that $\dom(F^{-1})=\rge(F)$ and that $(\ox, \oy)\in \ri(\gph(F))$ if and only if $(\oy, \ox)\in \ri(\gph(F^{-1}))$. Thus the conclusion follows directly from Theorem \ref{TheoRoc1}. $\h$

\section{Relative Interiors and Coderivatives of Generalized Epigraphical Mappings}
\setcounter{equation}{0}

Let us now apply Theorem \ref{TheoRoc1} to derive a representation for relative interiors of \emph{generalized epigraphical mappings} defined by
\begin{equation}\label{GEM}
F(x):=E_{f_1}(x)\times E_{f_2}(x)\times \cdots \times E_{f_m}(x), \; x\in \R^n,
\end{equation}
where $E_{f_i}(x):=[f_i(x), \infty)$ for $x\in \R^n$, and $f_i\colon \R^n\to (-\infty, \infty]$ for $i=1, \ldots, m$ are extended-real-valued convex functions.

\begin{Theorem}{\bf(relative interiors of generalized convex epigraphical graphs).}\label{TheoRoc2}
  Let $f_i\colon \R^n\to(-\infty,\infty]$ for $\ i=1,\ldots, m$ be extended-real-valued convex functions satisfying
  \begin{equation*}
    \bigcap_{i=1}^m\ri(\dom(f_i))\ne\emptyset.
  \end{equation*}
 Then we have the following representation for the generalized epigraphical mapping \eqref{GEM}:
  \begin{eqnarray*}
  \begin{array}{ll}
    \ri\big(\gph(F)\big)=\Big\{(x,\lambda_1,\ldots,\lambda_m)\in\R^n\times\R^m\; \big |\; &x\in \bigcap_{i=1}^{m}\ri(\dom(f_i)),\\
    &\;f_i(x)<\lambda_i\; \mbox{\rm for all } i=1,\ldots,m\Big\}.
    \end{array}
  \end{eqnarray*}
 \end{Theorem}
{\bf Proof.} It follows from the definition of the generalized epigraphical mapping \eqref{GEM} that $\dom(F)=\bigcap_{i=1}^m\dom(f_i)$ and
\begin{align*}
\gph(F)&=\big\{(x,\lambda)\in\R^n\times\R^m\; \big|\; x\in \dom(F), \lambda\in F(x)\big\}\\
 &=\big\{(x,\lambda_1,\ldots,\lambda_m)\in\R^n\times\R^m\; \big|\; x\in \bigcap_{i=1}^m\dom(f_i), \;f_i(x)\leq \lambda_i\; \mbox{\rm for all}\; i=1,\ldots, m\big\}.
 \end{align*}
For any $x\in\dom(F)$, it readily follows that
\begin{equation*}
  \ri(F(x))=\ri\big([f_1(x),\infty)\times\cdots\times[f_m(x),\infty)\big)=(f_1(x),\infty)\times\cdots\times(f_m(x),\infty).
\end{equation*}
Under the assumption that $ \bigcap_{i=1}^m\ri(\dom(f_i))\ne\emptyset$, we employ Theorem \ref{Ri_nonem} and get
\begin{equation*}
  \ri\Big( \bigcap_{i=1}^{m}\dom(f_i)\Big)= \bigcap_{i=1}^{m}\ri(\dom(f_i)).
\end{equation*}
 Applying finally Theorem~\ref{TheoRoc1} gives us
\begin{align*}
  \ri(\gph(F))&=\big\{(x,\lambda_1,\ldots,\lambda_m)\in\R^n\times\R^m\;\big |\; x\in \ri\big(\bigcap_{i=1}^{m}(\dom(f_i)\big), \; (\lambda_1,\ldots,\lambda_m)\in \ri(F(x))\big\}\\
    &=\big\{(x,\lambda_1,\ldots,\lambda_m)\in\R^n\times\R^m\; \big |\; x\in \bigcap_{i=1}^{m}\ri(\dom(f_i)), \;f_i(x)<\lambda_i\; \mbox{\rm for all } i=1,\ldots,m\big\},
\end{align*}
which thus completes the proof of this theorem. $\h$

As a direct consequence of Theorem \ref{TheoRoc1}, we obtain in the corollary below a representation for the relative interior of the epigraph of an extended-real-valued convex function; see, e.g., \cite[Proposition 1.1.9]{HU1}.

\begin{Corollary}{\bf (relative interiors of convex epigraphs).}\label{Co-Roc}
  Let $f\colon \R^n\to (-\infty, \infty]$ be an extended-real-valued convex function. Then
  \begin{equation*}\label{re-epi}
    \ri(\epi(f))=\big\{(x,\lambda)\in \R^n\times \R\;\big |\; x\in\ri(\dom(f)),\; f(x)<\lambda\big\}.
  \end{equation*}
\end{Corollary}

Given an extended-real-valued convex function $f\colon \R^n\to(-\infty,\infty]$, for each $\bar{x}\in\dom(f)$ we define the operation
\begin{equation*}
  \alpha\odot \partial f(\bar{x}):=
  \begin{cases}
  \alpha \partial f(\bar{x})&\; \text{if}\; \alpha >0,\\
  \partial^{\infty}f(\bar{x})&\; \text{if}\; \alpha=0,
  \end{cases}
\end{equation*}
where $\partial f(\bar{x})$ is the  subdifferential of $f$ at $\bar{x}$ defined in \eqref{cvsd}
and where $\partial^{\infty}f(\bar{x})$ is the {\em singular subdifferential} of $f$ at $\bar{x}$ defined by
  \begin{equation*}
    \partial^{\infty}f(\bar{x}):=\big\{v\in\R^n\;\big |\; (v,0)\in N((\bar{x},f(\bar{x}));\epi(f))\big\}.
  \end{equation*}
In the proof of the next result, we employ a well-known representation of the subdifferential of a convex function via the normal cone to its epigraph saying that
\begin{equation*}
  \partial f(\bar{x})=\big\{v\in\R^n\;\big |\; (v,-1)\in N((\bar{x},f(\bar{x}));\epi(f))\big\}, \; \ox\in \dom(f).
\end{equation*}

\begin{Proposition}\label{CEPG}{\bf (coderivative of epigraphs for extended-real-valued functions).}\label{lm-co-se-map}
  Let $f\colon\R^n\to(-\infty,\infty]$ be an extended-real-valued convex function, and let $F\colon\R^n\tto{\R}$ be the set-valued mapping defined by $F(x):=[f(x),\infty)$ for $x\in \R^n$. Then for any $\bar{x}\in\dom(F)=\dom(f)$ we have the representation
  \begin{equation*}\label{co-ex-set-map}
    D^*F(\bar{x},f(\bar{x}))(\alpha)=\left\{\begin{aligned}&\alpha\odot \partial f(\bar{x})&\ &\text{if}&\ &\alpha\geq 0,&\\
     &\emptyset&\ &\text{if}&\ &\alpha <0.&
    \end{aligned}\right.
  \end{equation*}
\end{Proposition}
{\bf Proof.} It is easy to see that $\gph(F)=\epi(f).$ The coderivative definition yields
\begin{equation}\label{eq-co-set}
  D^*F(\bar{x},f(\bar{x}))(\alpha)=\big\{v\in\R^n\;\big |\; (v,-\alpha)\in N((\bar{x},f(\bar{x}));\epi(f))\big\},\;\alpha\in\R.
\end{equation}
We consider the following three possible choices for $\alpha$:\\[1ex]
$\bullet$ If $\alpha>0,$ then $(v,-\alpha)\in N((\bar{x},f(\bar{x}));\epi(f))$ if and only if
\begin{equation*}
  \left(\frac{v}{\alpha},-1\right)\in N\big((\bar{x},f(\bar{x}));\epi(f)\big).
\end{equation*}
This means that $v/\alpha\in\partial f(\ox),$ and hence $v\in \alpha\partial f(\bar{x}).$\\[1ex]
$\bullet$ If $\alpha=0,$ then the definition of the singular subdifferential yields $(v,-\alpha)\in N((\bar{x},f(\bar{x}));\epi(f))$ if and only if $v\in\partial^{\infty}f(\ox)$.\\[1ex]
$\bullet$ If $\alpha<0$, then we can show that $D^*F(\bar{x},f(\bar{x}))(\alpha)=\emptyset$. Indeed, suppose that this is not the case and find $(v,-\alpha)\in N((\bar{x},f(\bar{x}));\epi(f))$.  Then
\begin{equation*}
\la v, x-\ox\ra-\alpha\big(\lambda-f(\ox)\big)\leq 0\; \mbox{\rm whenever }(x, \lambda)\in \epi(f).
\end{equation*}
Choosing $(\bar{x},f(\bar{x})+1)\in\epi(f)$, we can see that $-\alpha\leq 0$. This contradiction verifies that $D^*F(\bar{x},f(\bar{x}))(\alpha)=\emptyset$.

Therefore, the representation of $D^*F(\bar{x},f(\bar{x}))(\alpha)$ follows from \eqref{eq-co-set} and the above definition of the operation $\alpha\odot \partial f(\bar{x})$.
$\h$

Now we are ready to obtain a useful coderivative representation for generalized epigraphical set-valued mappings under the relative interior condition.

\begin{Theorem}\label{CGEM}{\bf(coderivatives of generalized epigraphical mappings).}
Let $f_i\colon \R^n\to(-\infty,\infty]$ for  $i=1,\ldots,m$ be extended-real-valued convex functions. Consider the generalized epigraphical set-valued mapping $F$ defined in \eqref{GEM}  and suppose that
\begin{equation*}
  \bigcap_{i=1}^m\ri(\dom(f_i))\neq\emptyset.
  \end{equation*}
   Then for any $\bar{x}\in \dom(F)$  we have the representation
\begin{equation*}
  D^*F(\bar{x},\bar{y})(\alpha)=\left\{
  \begin{aligned} &\sum_{i=1}^{m}\alpha_i\odot \partial f_i(\bar{x})&\ &\text{if}&\ &\alpha_i\geq0\; \mbox{\rm for all }\ i=1,\ldots m,&\\
  &\emptyset&&\text{if}& &\alpha_i<0\; \mbox{\rm for some }i=1, \ldots, m,&
  \end{aligned}\
  \right.
\end{equation*}
where $\bar{y}:=(f_1(\ox), \ldots, f_m(\ox))$ and $\alpha:=(\alpha_1,\ldots,\alpha_m)\in \R^m.$
\end{Theorem}
{\bf Proof.}
Define the sets
\begin{equation*}
  \Omega_i:=\big\{(x,\lambda_1,\ldots,\lambda_m)\in \R^n\times \R^m\;\big |\; x\in \R^n, \; \lambda_i\geq f_i(x)\big\},\; i=1,\ldots,m,
\end{equation*}
and observe that $\gph(F)=\bigcap^m_{i=1}\Omega_i$. It follows from Corollary~\ref{Co-Roc}  that
\begin{equation*}
  \ri(\Omega_i)=\big\{(x,\lambda_1,\ldots,\lambda_m)\in \R^n\times \R^m\; \big|\; x\in \ri(\dom(f_i)), \lambda_i> f_i(x)\big\},\; i=1,\ldots,m.
\end{equation*}
Choose $x_0\in \bigcap_{i=1}^m\ri(\dom(f_i))$ and let $\bar{\lambda}_i:=f_i(x_0)+1>f_i(x_0)$ for $i=1, \ldots, m$. Then $(x_0, \bar{\lambda}_1, \ldots, \bar{\lambda}_m)\in \bigcap_{i=1}^{m}\ri(\Omega_i)$, and hence $\bigcap_{i=1}^{m}\ri(\Omega_i)\ne\emptyset.$

Applying now the normal cone intersection rule from Theorem~\ref{ncr} gives us
\begin{equation*}  N\big((\bar{x},\bar{y});\gph(F))\big)=N\big((\bar{x},\bar{y});\bigcap_{i=1}^{m}\Omega_i)\big) =\sum_{i=1}^{m}N\big((\bar{x},\bar{y});\Omega_i\big).
\end{equation*}
It is easy to see that $\Omega_1=\epi(f_1)\times\R^{m-1}$, which gives us by \cite[Proposition 2.11]{bmn} that
\begin{equation*}
N\big((\bar{x},\bar{y});\Omega_1\big)=N\big((\bar{x},f_1(\bar{x}));\epi(f_1)\big)\times \{0\}.
\end{equation*}
The latter means that  $(v,-\alpha)\in N\big((\bar{x},\bar{y});\Omega_1\big)$ if and only if $(v, -\alpha_1)\in N\big((\bar{x},f_1(\bar{x}));\epi(f_1)\big)$ and $\alpha_j=0$ for $j=2, \ldots, m$. In general, we observe that $(v,-\alpha)\in N\big((\bar{x},\bar{y});\Omega_i\big)$
if and only if $(v, -\alpha_i)\in N\big((\bar{x},f_i(\bar{x}));\epi(f_i)\big)$ and $\alpha_j=0$ for $j\in \{1, \ldots, m\}\setminus\{i\}$.

Finally, it follows from the coderivative construction and the proof of Proposition~\ref{lm-co-se-map} that
\begin{equation*}
\begin{aligned}
  D^*F(\bar{x},\bar{y})(\alpha)&=\Big\{v\in \R^n\; \Big|\; (v,-\alpha)\in N((\bar{x},\bar{y}); \gph(F))\Big\}\\
  &=\Big\{v\in \R^n\ |\ (v,-\alpha)\in \sum_{i=1}^{m}N\big((\bar{x},\bar{y});\Omega_i)\big)\Big\}\\
&=\left\{
  \begin{aligned} &\sum_{i=1}^{m}\alpha_i\odot \partial f_i(\bar{x})&\ &\text{if}&\ &\alpha_i\geq0\ \text{for all}\ i=1,\ldots ,m,&\\
  &\emptyset&&\text{if}& &\alpha_i<0\; \mbox{\rm for some }i=1, \ldots, m,&
  \end{aligned}\
  \right.
  \end{aligned}
\end{equation*}
which therefore completes the proof of the theorem. $\h$

\section{Optimal Value Functions and Generalized Chain Rules}
\setcounter{equation}{0}

Given a  set-valued mapping $F\colon \R^n\tto \R^m$ and an extended-real-valued function $\ph\colon \R^m\to (-\infty, \infty]$, define the associated \emph{optimal value function} by
\begin{equation}\label{mgn}
\mu(x):=\inf\big\{\ph(y)\; \big|\; y\in F(x)\big\},\; x\in \R^n.
\end{equation}
Throughout this section, we assume that $\mu(x)>-\infty$ for all $x\in \R^n$. For any $\ox\in \dom(\mu)$, consider the argminimum set
\begin{equation*}
S(\ox):=\big\{y\in F(\ox)\; \big|\; \mu(\ox)=\ph(y)\big\}.
\end{equation*}

We have the following exact formula for subdifferentiation of the optimal value function under the relative interior qualification condition.

\begin{Proposition}{\bf (subdifferentials of optimal value functions).}\label{smgn} Let $\mu$ be the optimal value function defined in \eqref{mgn}, where $F$ is a convex set-valued mapping, and where $\ph$ is an extended-real-valued convex function. Then the function $\mu$ is convex. In addition, for any $\ox\in \dom(\mu)$ and any $\oy\in S(\ox)$ we have
\begin{equation*}
\partial \mu(\ox)=\bigcup_{v\in \partial \ph(\oy)}D^*F(\ox, \oy)(v)
\end{equation*}
provided that there exists $x_0\in \ri(\dom(F))$ such that $\ri(F(x_0))\cap \ri(\dom(\ph))\neq \emptyset$.
\end{Proposition}
{\bf Proof.} Define the function $\psi\colon \R^n\times \R^m\to (-\infty, \infty]$ by $\psi(x,y):=\ph(y)$ for $(x, y)\in \R^n\times \R^m$. Then $\psi$ is clearly convex with $\dom(\psi)=\R^n\times \dom(\ph)$, and hence we get $\ri(\dom(\psi))=\R^n\times \ri(\dom(\ph))$. Choose  $y_0\in \ri(F(x_0))\cap \ri(\dom(\ph))$. Then $(x_0, y_0)\in \ri(\dom(\psi))$, and it follows from  Theorem~\ref{TheoRoc1} that $(x_0, y_0)\in \ri(\gph(F))$. Therefore, $\ri(\gph(F))\cap \ri(\dom(\psi))\neq \emptyset$. Now we deduce the claimed result from \cite[Theorem 9.1]{bmncal}. $\h$

Recall that $\ph\colon \R^m\to (-\infty, \infty]$ is \emph{nondecreasing componentwise} if we have
\begin{equation*}
\big[x_i\leq u_i\; \mbox{\rm for all }i=1, \ldots, m]\Longrightarrow [\ph(x_1, \ldots, x_m)\leq \ph(u_1, \ldots, u_m)\big].
\end{equation*}
The next theorem gives us a generalization of \cite[Theorem 4.3.1]{HU1} for a broad class of composite extended-real-valued functions. We also provide a simpler proof of this theorem applying the coderivative of generalized epigraphical mappings.

\begin{Theorem}{\bf(subdifferentials of a composition with increasing extended-real-valued functions of several variables).} Let $f_i\colon \R^n\to\R$ for $i=1, \ldots, m$ be real-valued convex functions, and let $\ph\colon \R^m\to (-\infty, \infty]$ be nondecreasing componentwise and convex. Consider the composite function
\begin{equation*}
g(x):=\ph(f_1(x), \ldots, f_m(x)), \; x\in \R^n.
\end{equation*}
Then the function $g\colon \R^n\to (-\infty, \infty]$ is convex. Suppose in addition  that there exists $(x_0, \lambda_1, \ldots, \lambda_m)\in \R^n\times \R^m$ such that $\lambda_i>f_i(x_0)$ for all $i=1, \ldots, m$ and $(\lambda_1, \ldots, \lambda_m)\in \ri(\dom(\ph))$. Then for any $\ox\in \dom(g)$ we have the subdifferential formula
\begin{equation*}
\partial g(\ox)=\bigg\{\sum_{i=1}^m \gamma_i \partial f_i(\ox)\;\bigg|\; (\gamma_1, \ldots, \gamma_m)\in \partial \ph(\oy)\bigg\},
\end{equation*}
where $\oy:=(f_1(\ox), \ldots, f_m(\ox))$.
\end{Theorem}
{\bf Proof.} Define the set-valued mapping $F(x):=[f_1(x), \infty)\times \cdots\times [f_m(x), \infty)$ for $x\in \R^n$ and then deduce from the nondecreasing componentwise property of $\ph$ that
\begin{equation*}
\mu(x)=g(x)\; \mbox{\rm for all }x\in \R^n,
\end{equation*}
where $\mu$ is the optimal value function \eqref{mgn} generated by $F$ and $\ph$. Observe that in this case we have that each function $f_i$ is continuous, and that $\gamma\odot \partial f_i(\ox)=\gamma \partial f(\ox)$ whenever $\gamma\geq 0$ and $\ox\in \R^n$. Using the representation
\begin{equation*}
\ri(F(x))=(f_1(x), \infty) \times \cdots\times (f_m(x), \infty)\;\mbox{ for any }\;x\in\R^n,
\end{equation*}
it follows from the imposed assumptions that there exists $x_0\in \ri(\dom(F))=\R^n$ such that $\ri(F(x_0))\cap \ri(\dom(\ph))\neq \emptyset$. Furthermore, Proposition~\ref{smgn} tells us that
\begin{equation*}
\partial g(\bar{x})=\partial \mu(\ox)=\bigcup_{\gamma\in \partial \ph(\oy)}D^*F(\ox, \oy)(\gamma).
\end{equation*}
The rest of the proof follows from the coderivative formula for $F$ in Theorem \ref{CGEM}. $\h$

\section{Coderivative Calculus in Finite-Dimensional Spaces}
\setcounter{equation}{0}

In this section, under the relative interior conditions imposed on domains and ranges of mappings, we establish major formulas of coderivative calculus including  sum rule, chain rule, and intersection rule for set-valued mappings in finite-dimensional spaces. The obtained results improve those  in \cite{bmncal} under more restrictive qualification conditions.

Given two set-valued mappings $F_1,F_2\colon\R^n\tto\R^m$, their {\em sum} is defined by
\begin{equation*}
(F_1+F_2)(x)=F_1(x)+F_2(x):=\big\{y_1+y_2\;\big|\;y_1\in F_1(x),\;y_2\in F_2(x)\big\}.
\end{equation*}
It is easy to see that $\dom(F_1+F_2)=\dom(F_1)\cap\dom(F_2)$, and that $F_1+F_2$ is convex provided that both $F_1$ and $F_2$ have this property.

Our first calculus result concerns representing coderivatives of sums $F_1+F_2$ at a given point $(\ox,\oy)\in\gph(F_1+F_2)$. To formulate this result, consider the nonempty set
\begin{equation*}
S(\ox,\oy):=\big\{(\oy_1,\oy_2)\in\R^m\times\R^m\;\big|\;\oy=\oy_1+\oy_2,\;\oy_i\in F_i(\ox)\;\mbox{\rm for }\;i=1,2\big\}.
\end{equation*}
The following theorem gives us the coderivative sum rule for set-valued mappings on finite-dimensional spaces. In this version, we use the relative interior qualification condition on domains replacing the condition on graphs known from \cite[Theorem~11.1]{bmncal}.

\begin{Theorem}{\bf(coderivative sum rule via qualification condition on domains).}\label{sumriule} Let $F_1,F_2\colon\R^n\tto\R^m$ be convex set-valued mappings. Imposing the relative interior condition
\begin{equation}\label{QC}
\ri(\dom(F_1))\cap\ri(\dom(F_2))\ne\emptyset,
\end{equation}
we have the coderivative sum rule
\begin{equation*}\label{csr}
D^*(F_1+F_2)(\ox,\oy)(v)=\bigcap_{(\oy_1,\oy_2)\in S(\ox,\oy)}\big[D^*F_1(\ox,\oy_1)(v)+D^*F_2(\ox,\oy_2)(v)\big]
\end{equation*}
for all $(\ox,\oy)\in\gph(F_1+F_2)$ and $v\in\R^m$.
\end{Theorem}
{\bf Proof.} Fix any $u\in D^*(F_1+F_2)(\ox,\oy)(v)$ and $(\oy_1,\oy_2)\in S(\ox,\oy)$ for which we have the inclusion $(u,-v)\in N((\ox,\oy);\gph(F_1+F_2))$.
Consider the convex sets
\begin{align*}
&\Omega_1:=\big\{(x,y_1,y_2)\in \R^n\times\R^m\times\R^m\;\big|\;y_1\in F_1(x)\big\},\\
&\Omega_2:=\big\{(x,y_1,y_2)\in\R^n\times\R^m\times\R^m\;\big|\;y_2\in F_2(x)\big\}
\end{align*}
and deduce from the normal cone definition that
\begin{equation*}
(u,-v,-v)\in N((\ox,\oy_1,\oy_2);\Omega_1\cap\Omega_2).
\end{equation*}
Now we intend to verify the inclusion
\begin{equation}\label{Eq-Sum-Rule-Nor}
(u,-v,-v)\in N((\ox,\oy_1,\oy_2);\Omega_1)+N((\ox,\oy_1,\oy_2);\Omega_2).
\end{equation}
Indeed, it follows from \eqref{QC} that there exists $x\in \ri(\dom(F_1))\cap\ri(\dom(F_2)),$ and hence Theorem \ref{Ri_nonem} implies that $\ri(F_1(x))\ne\emptyset$ and $\ri(F_2(x))\ne\emptyset$. Theorem \ref{TheoRoc1} ensures that
\begin{align*}
&\ri(\Omega_1)=\big\{(x,y_1,y_2)\in\R^n\times\R^m\times\R^m\;\big|\;x\in\ri(\dom (F_1)),\; y_1\in\ri(F_1(x))\big\},\\
&\ri(\Omega_2)=\big\{(x,y_1,y_2)\in\R^n\times\R^m\times\R^m\;\big|\;x\in\ri(\dom (F_2)), \;y_2\in\ri(F_2(x))\big\},
\end{align*}
which shows in turn that condition \eqref{QC} yields $\ri(\Omega_1)\cap\ri(\Omega_2)\ne\emp$. This tells us by Theorem~\ref{ncr} that \eqref{Eq-Sum-Rule-Nor} is satisfied, and therefore we get the relationships
\begin{equation*}
(u,-v,-v)=(u_1,-v,0)+(u_2,0,-v)\textrm{ with }(u_i,-v)\in N((\ox,\oy_i);\gph(F_i))\textrm{ for }i=1,2.
\end{equation*}
This implies by the coderivative definition that
\begin{equation*}
u=u_1+u_2\in D^*F_1(\ox,\oy_1)(v)+D^*F_2(\ox,\oy_2)(v)
\end{equation*}
as desired. The reverse inclusion is obvious, and thus we verify the claimed sum rule. $\h$\vspace*{0.05in}

Next we present the well-known subdifferential sum rule (see, e.g., \cite[Theorem~23.8]{r}), which can also be derived from Theorem~\ref{sumriule}.

\begin{Corollary} {\bf(subdifferential sum rule).}\label{sr}
Let $f_i\colon \R^n\to\oR$, $i=1,2$, be  extended-real-valued convex functions. Suppose that the relative interior qualification condition
\begin{equation}\label{LQC}
\ri(\dom(f_1))\cap\ri(\dom(f_2))\ne\emptyset
\end{equation}
is satisfied.  Then for all $\ox\in\dom(f_1)\cap\dom(f_2)$ we have the subdifferential sum rule
\begin{equation}\label{ssr}
\partial(f_1+f_2)(\ox)=\partial f_1(\ox)+\partial f_2(\ox).
\end{equation}
\end{Corollary}
{\bf Proof.} Define the convex set-valued mappings $F_1,F_2\colon X\tto\R$ by
\begin{equation*}
F_i(x):=\big[f_i(x),\infty\big)\;\mbox{ for }\;i=1,2.
\end{equation*}
It is easy to see that $\gph(F_i)=\epi(f_i)$ and $\dom(F_i)=\dom(f_i)$ for $i=1,2$.  Furthermore, the qualification condition \eqref{LQC} clearly implies the fulfillment of \eqref{QC}.

To proceed further, fix any $\ox\in\dom(f_1)\cap\dom(f_2)$, and let $\oy:=f_1(\ox)+f_2(\ox)$. For every $x^*\in\partial(f_1+f_2)(\ox)$ we have the coderivative inclusion
\begin{equation*}
x^*\in D^*(F_1+F_2)(\ox,\oy)(1).
\end{equation*}
Applying to the latter Theorem~\ref{sumriule} with $\oy_i=f_i(\ox)$ for $i=1,2$ gives us
\begin{equation*}
x^*\in D^*F_1(\ox,\oy_1)(1)+D^*F_2(\ox,\oy_2)(1)=\partial f_1(\ox)+\partial f_2(\ox),
\end{equation*}
which verifies the inclusion ``$\subset$" in \eqref{ssr}. The reverse inclusion is obvious.
$\h$\vspace*{0.05in}

Now we define the {\em composition} of two set-valued mappings $F\colon\R^n\tto\R^m$ and $G\colon\R^m\tto\R^q$ by
\begin{equation*}
(G\circ F)(x)=\bigcup_{y\in F(x)}G(y):=\big\{z\in G(y)\;\big|\;y\in F(x)\big\},\; x\in\R^n,
\end{equation*}
and observe that $G\circ F$ is convex provided that both $F$ and $G$ have this property. Given $\oz\in(G\circ F)(\ox)$, we consider the set
\begin{equation*}
M(\ox,\oz):=F(\ox)\cap G^{-1}(\oz).
\end{equation*}
The following theorem establishes the coderivative chain rule for set-valued mappings in finite-dimensional spaces. In this version, we use the relative interior qualification condition on domains and ranges replacing the one on graphs known from \cite[Theorem~11.2]{bmncal}.

\begin{Theorem}{\bf(coderivative chain rule via qualification condition on domains).}\label{scr} Let $F\colon\R^n\tto\R^m$ and $G\colon\R^m\tto\R^q$ be convex set-valued mappings satisfying the relative interior qualification
condition
\begin{equation}\label{QCCR}\ri\big(\rge(F)\big)\cap
\ri\big(\dom(G)\big)\ne\emp.
\end{equation}
Then for any $(\ox,\oz)\in\gph(G\circ F)$ and $w\in\R^q$ we have the coderivative chain rule
\begin{equation}\label{chain}
D^*(G\circ F)(\ox,\oz)(w)=\bigcap_{\oy\in M(\ox,\oz)}D^*F(\ox,\oy)\circ D^*G(\oy,\oz)(w).
\end{equation}
\end{Theorem}
{\bf Proof.} Picking $u\in D^*(G\circ F)(\ox,\oz)(w)$ and $\oy\in M(\ox,\oz)$ gives us the inclusion $(u,-w)\in N((\ox,\oz);\gph(G\circ F))$, which means that
\begin{equation*}
\la u,x-\ox\ra-\la w,z-\oz\ra\le 0\;\mbox{\rm for all }\;(x,z)\in\gph(G\circ F).
\end{equation*}
Define two convex subsets of $\R^n\times\R^m\times\R^q$ by
$$
\Omega_1:=\gph(F)\times\R^q\;\mbox{ and }\;\Omega_2:=\R^n\times\gph(G).
$$
It is easy to see that
\begin{equation}\label{graphrelation1}
\Omega_1-\Omega_2=\R^n\times
\big(\rge(F)-\dom(G)\big)\times \R^q.
\end{equation}
Using \eqref{graphrelation1}, we have the representation \begin{equation*}\ri(\Omega_1-\Omega_2)=\R^n\times
\ri\big(\rge(F)-\dom(G)\big)\times \R^q. \end{equation*} It follows from \eqref{QCCR} due to the definitions of the sets $\Omega_1$ and $\Omega_2$ that
\begin{equation}\label{QC1_re} 0\in\ri(\Omega_1-\Omega_2),\;\mbox{and
so }\;\ri(\Omega_1)\cap\ri(\Omega_2)\ne\emp.\end{equation}

We can directly deduce from the definitions that
\begin{equation*}
(u,0,-w)\in N((\ox,\oy,\oz);\Omega_1\cap\Omega_2).
\end{equation*}
Applying Theorem~\ref{ncr} with qualification \eqref{QC1_re} tells us that
\begin{equation*}
(u,0,-w)\in N((\ox,\oy,\oz);\Omega_1\cap\Omega_2)=N((\ox,\oy,\oz);\Omega_1)+N((\ox,\oy,\oz);\Omega_2).
\end{equation*}
Thus using further the definitions of the sets $\Omega_1$ and $\Omega_2$ based on $\gph(F)$ and $\gph(G)$, there exists a vector $v\in\R^m$ such that we have the representation
\begin{equation*}
(u,0,-w)=(u, -v,0)+(0,v,-w),
\end{equation*}
where $(u,-v)\in N((\ox,\oy);\gph(F)),\;(v,-w)\in N((\oy,\oz);\gph(G))$.
This shows by the coderivative definition (\ref{cod}) that
\begin{equation*}
u\in D^*F(\ox,\oy)(v)\;\mbox{\rm and}\;v\in D^*G(\oy,\oz)(w),
\end{equation*}
and so we verify the inclusion ``$\subset$" in \eqref{chain}. The reverse inclusion is trivial. $\h$\vspace*{0.05in}

Let $F\colon\R^n\tto \R^m$ be a set-valued mapping, and let $\Theta\subset \R^m$ be a given set. The {\em preimage} or {\em inverse image} of $\Theta$ under the mapping $F$ is given by
\begin{equation*}
F^{-1}(\Theta)=\big\{x\in \R^n\;\big |\; F(x)\cap \Theta\ne\emptyset\big\}.
\end{equation*}
The next result gives us a representation of the normal cone to $F^{-1}(\Theta)$ via the normal cone of $\Theta$ and the coderivative of $F$. We use here the relative interior qualification condition on ranges replacing the condition on graphs known from \cite[Proposition~10.1]{bmncal}.

\begin{Proposition}\label{IVIM}{\bf(representation of the normal cone to preimages).}\label{Theo-code-preimage}
Let $F\colon\R^n\tto\R^m$ be a convex set-valued mapping, and let $\Theta\subset\R^m$ be a convex set. Suppose that
\begin{equation}\label{IRQC}
\ri\big(\rge(F)\big)\cap\ri(\Theta)\ne\emptyset.
\end{equation}
Then for any $\bar{x}\in F^{-1}(\Theta)$ and $\bar{y}\in F(\bar{x})\cap \Theta$ we have the representation
\begin{equation*}
N(\bar{x};F^{-1}(\Theta))=D^*F(\bar{x},\bar{y})\big(N(\bar{y};\Theta)\big).
\end{equation*}
\end{Proposition}
{\bf Proof.}
In the setting of Theorem~\ref{scr}, consider the set-valued mapping $G\colon \R^m\tto \R^q$ given by
\begin{equation*}
G(x):=\Delta_{\Theta}(x) =
\begin{cases}
0 & \mbox{\rm if}\;  x\in\Theta, \\
\emptyset &\mbox{\rm if}\; x\notin\Theta.
\end{cases}
\end{equation*}
It is clear that $\dom(\Delta_\Theta)=\Theta$, $\gph(\Delta_{\Theta})=\Theta\times\{0\}$ and that for any $\ox\in\Theta$ we get $N\big((\ox,0);\gph(\Delta_{\Theta})\big)=N(\ox;\Theta)\times \R^q$. Therefore,
\begin{equation*}
D^*\Delta_{\Theta}(\ox,0)(v)=N(\ox;\Theta)\; \text{for all}\; v\in \R^q.
\end{equation*}
It is easy to check the composite representation
\begin{equation*}
\Delta_{F^{-1}(\Theta)}(x)=\big(\Delta_{\Theta}\circ F\big)(x)\; \text{for all}\; x\in\R^n,
\end{equation*}
where $\Delta_{F^{-1}(\Theta)}\colon \R^n\tto\R^p$ is given by $\Delta_{F^{-1}(\Theta)}(x):=0$ if $x\in F^{-1}(\Theta)$, and $\Delta_{F^{-1}(\Theta)}(x):=\emptyset$ otherwise. Observe that the imposed relative interior qualification condition \eqref{IRQC} guarantees that $\ri(\rge(F))\cap \ri(\dom(G))\neq\emptyset$.  Then the claimed formula for $N(\ox;F^{-1}(\Theta))$ follows from the coderivative chain rule of Theorem~\ref{scr} with the outer mapping $G:=\Delta_{\Theta}.$
$\h$\vspace*{0.05in}

The last result of this section provides a precise representation formula for the normal cone to sublevel sets of extended-real-valued convex functions.

\begin{Corollary}{\bf(representation of the normal cone to sublevel sets).}
Let $f\colon\R^n\to\oR$ be a convex function. For $\lambda\in\R$, consider  the sublevel set
\begin{equation*}
\mathcal{L}_\lambda:=\big\{x\in\R^n\;\big|\;f(x)\le\lambda\big\}
\end{equation*}
Assume that $f(\ox)=\lambda$, and that there exists $\hat{x}\in \ri(\dom(f))$ such that $f(\hat{x})<\lambda$. Then we have the representation
\begin{equation*}
N(\ox;\mathcal{L}_\lambda)=\disp\bigcup_{\alpha\ge 0}\alpha\odot\partial f(\ox).
\end{equation*}
\end{Corollary}
{\bf Proof.} Let $F(x):=E_f(x)$ for $x\in \R^n$, and let $\Theta:=(-\infty, \lambda]$. Then $\mathcal{L}_\lambda=F^{-1}(\Theta)$. Since $\hat{x}\in \ri(\dom(f))$ and  $f(\hat{x})<f(\ox)=\lambda$, we can show that
$$\ri\big(\rge(F)\big)\cap\ri(\Theta)\ne\emptyset.$$
Indeed, choose $\gamma\in\R$ such that $f(\ox)<\gamma<\lambda$. By Corollary \ref{Co-Roc}, we see that $(\ox, \gamma)\in \ri(\epi(f))=\ri(\gph(F))$, so $\gamma\in \ri(\rge(F))$ by Corollary \ref{rirange}. Thus  $\gamma\in \ri\big(\rge(F)\big)\cap\ri(\Theta)$. Since $N(f(\ox); \Theta)=N(\lambda; \Theta)=[0, \infty)$, by Proposition \ref{CEPG} and Proposition~\ref{IVIM} we have
\begin{eqnarray*}
\begin{array}{ll}
N(\ox;\mathcal{L}_\lambda)&=N(\bar{x};F^{-1}(\Theta))=D^*F(\bar{x},\bar{y})\big(N(\lambda;\Theta)\big)\\
&=D^*F(\ox, f(\ox))\big([0, \infty)\big)=\disp\bigcup_{\alpha\ge 0}D^*F(\ox, f(\ox))(\alpha)=\disp\bigcup_{\alpha\ge 0}\alpha\odot\partial f(\ox),
\end{array}
\end{eqnarray*}
which completes the proof of the corollary. $\h$\\[1ex]
{\bf Acknowledgements}. The authors are very grateful to
the anonymous referees for their valuable remarks and suggestions
that allowed us to improve the original presentation.

\end{document}